%&amstex
\input amstex
\input amsppt.sty
\magnification=\magstep1
\hsize=30truecc
\baselineskip=16truept
\vsize=22.2truecm
\NoBlackBoxes
\nologo
\pageno=1
\topmatter
\TagsOnRight

\def\N{\Bbb N}
\def\Z{\Bbb Z}

\def\C{\Bbb C}
\def\l{\left}
\def\r{\right}
\def\b{\bigg}
\def\bg{\bigg}
\def\({\b(}
\def\[{\b[}
\def\){\b)}
\def\]{\b]}

\def\t{\text}
\def\f{\frac}
\def\mo{\roman{mod}}

\def\se {\subseteq}

\def\sm{\setminus}

\def\bi{\binom}
\def\eq{\equiv}
\def\cs{\cdots}
\def\ls{\leqslant}
\def\gs{\geqslant}

\def\ve{\varepsilon}

\def\Proof{\noindent{\it Proof}}
\def\Remark{\noindent{\it Remark}}

\hbox{J. Combin. Theory Ser. A 116(2009), no.\,8, 1374--1381.}
\medskip
\title A new extension of the Erd\H os-Heilbronn conjecture\endtitle
\author Hao Pan and Zhi-Wei Sun \endauthor

\affil Department of Mathematics, Nanjing University
     \\Nanjing 210093, People's Republic of China
    \\ {\tt haopan79\@yahoo.com.cn},  \quad {\tt zwsun\@nju.edu.cn}
 \endaffil
\abstract Let $A_1, \ldots, A_n$ be finite subsets
of a field $F$, and let
$$f(x_1,\ldots,x_n)=x_1^k+\cdots+x_n^k+g(x_1,\ldots,x_n)\in F[x_1,\ldots,x_n]$$
with $\deg g<k$. We obtain a lower bound for the cardinality of
$$\{f(x_1,\ldots,x_n):\, x_1\in A_1,\ldots,x_n\in A_n,\ \t{and}\ x_i\not=x_j\text{ if}\ i\not=j\}.$$
The result extends the Erd\H os-Heilbronn conjecture in a new way.
\endabstract
\thanks 2000 {\it Mathematics Subject Classification}.
Primary 11B75; Secondary 05A05, 11P99, 11T06, 12E10.
\newline\indent
The second author is responsible for communications, and supported
by the National Natural Science Foundation (grant 10871087) of
China.
\endthanks
\endtopmatter
\document

\heading{1. Introduction}\endheading

In 1964 P. Erd\H os and H. Heilbronn [EH] made the following
challenging conjecture: If $p$ is a prime, then for any
subset $A$ of the finite field $\Z/p\Z$ we have
$$|\{x_1+x_2:\, x_1,x_2\in A\ \t{and}\ x_1\not=x_2\}|\gs\min\{p,\,2|A|-3\}.$$
It had remained open for thirty years until it was confirmed fully
by Dias da Silva and Hamidoune [DH] who actually obtained the
following generalization with the help of the representation theory
of symmetric groups: For any finite subset $A$ of a field $F$, we
have the inequality
$$\align&|\{x_1+\cdots+x_n:\, x_1,\ldots,x_n\in A,\ \t{and}\ x_i\not=x_j\ \t{if}\ i\not=j\}|
\\&\qquad\qquad\quad\gs\min\{p(F),\,n(|A|-n)+1\},
\endalign$$
where $p(F)=p$ if the characteristic of $F$ is a prime $p$, and
$p(F)=+\infty$ if $F$ is of characteristic zero. Recently
P. Balister and J. P. Wheeler [BW] extended the Erd\H os-Heilbronn conjecture
to any finite group; namely they showed that for any nonempty
subsets $A_1$ and $A_2$ of a finite group $G$ written additively we have
$$|\{x_1+x_2:\, x_1\in A_1,\, x_2\in A_2\ \t{and}\ x_1\not=x_2\}|\gs\min\{p(G),\,|A_1|+|A_2|-3\},$$
where $p(G)$ is the least positive order of a nonzero element of $G$, and
$p(G)$ is regarded as $+\infty$ if $G$ is torsion-free.

In 1996 N. Alon, M. B. Nathanson and I. Z. Ruzsa [ANR] used the
so-called polynomial method (see also Alon [A], Nathanson [N, pp.\, 98--107],
and T. Tao and V. H. Vu [TV, pp.\,329--345]) to deduce the following result: If
$A_1,\ldots,A_n$ are finite subsets of a field $F$ with
$0<|A_1|<\cdots<|A_n|$, then
$$\align&|\{x_1+\cdots+x_n:\, x_1\in A_1,\ldots,x_n\in A_n,\ \t{and}\ x_i\not=x_j\ \t{if}\ i\not=j\}|
\\&\qquad\qquad\quad\gs\min\bg\{p(F),\,1+\sum_{i=1}^n(|A_i|-i)\bg\}.
\endalign$$
Consequently, if $A_1,\ldots,A_n$ are finite subsets of a field $F$
with $|A_i|\gs i$ for all $i=1,\ldots,n$, then
$$\align&|\{x_1+\cdots+x_n:\, x_1\in A_1,\ldots,x_n\in A_n,\ \t{and}\ x_i\not=x_j\ \t{if}\ i\not=j\}|
\\&\qquad\qquad\quad\gs\min\bg\{p(F),\,1+\sum_{i=1}^n\min_{i\ls j\ls n}(|A_j|-j)\bg\}.
\endalign$$
(Choose $A_i'\se A_i$ with $|A_i'|=i+\min_{i\ls j\ls n}(|A_j|-j)\ls
|A_i|$. Then $|A_1'|<\cdots<|A_n'|$.) For other results on
restricted sumsets obtained by the polynomial method, the reader may
consult [HS], [PS], [S03], [SY] and [S08a].

Recently Z. W. Sun [S08] obtained the following result on value sets of polynomials.

\proclaim{Theorem 1.1 {\rm (Z. W. Sun [S08])}} Let
$A_1,\ldots,A_n$ be finite nonempty subsets of a field $F$, and
let
$$f(x_1,\ldots,x_n)=a_1x_1^k+\cdots+a_nx_n^k+g(x_1,\ldots,x_n)\in F[x_1,\ldots,x_n]\tag 1.1$$
 with
 $$k\in\Z^+=\{1,2,3,\ldots\},\ a_1,\ldots,a_n\in F^*=F\sm\{0\}\ \t{and}\ \deg g<k.\tag1.2$$

 {\rm (i)} We have
 $$\aligned&|\{f(x_1,\ldots,x_n):\,x_1\in A_1,\ldots,x_n\in A_n\}|
 \\&\quad\gs\min\bg\{p(F),\,1+\sum_{i=1}^n\l\lfloor\f{|A_i|-1}k\r\rfloor\bg\}.
 \endaligned\tag1.3$$

 {\rm (ii)} If $k\gs n$ and $|A_i|\gs i$ for $i=1,\ldots,n$, then
 $$\aligned &|\{f(x_1,\ldots,x_n):\,x_1\in A_1,\ldots,x_n\in A_n,\ \t{and}\ x_i\not=x_j\ \t{if}\ i\not=j\}|
 \\&\qquad\qquad\gs\min\bg\{p(F),\,1+\sum_{i=1}^n\l\lfloor\f{|A_i|-i}k\r\rfloor\bg\}.
 \endaligned\tag1.4$$
\endproclaim

Throughout this paper, for a predicate $P$ we let
$$[\![P]\!]=\cases1&\t{if}\ P\ \t{holds},\\0&\t{otherwise}.\endcases$$
For $a\in\Z$ and $k\in\Z^+$, we use $\{a\}_k$ to denote the least nonnegative residue of $a$ modulo $k$.

 Let $\C$ be the field of complex numbers.
By [S08, Example 4.1], if $k\in\Z^+$, $q\in\{0,1,\ldots\}$, and
$A=\{z\in\C:\, z^k\in\{1,\ldots,q\}\}\cup R$ with $R\se\{z\in\C:\,
z^k=q+1\}$ and $|R|=r<k$, then $|A|=kq+r$ and
$$\align&|\{x_1^k+\cdots+x_n^k:\, x_1,\ldots,x_n\in A,\ \t{and}\ x_i\not=x_j\ \t{if}\ i\not=j\}|
\\&=\f{n(|A|-n)-\{n\}_k\{|A|-n\}_k}k+r[\![\{n\}_k>r]\!]+1.
\endalign$$
Motivated by this example, Sun [S08] raised the following extension of the Erd\H os-Heilbronn conjecture.

\proclaim{Conjecture 1.1 {\rm (Z. W. Sun [S08])}} Let
$f(x_1,\ldots,x_n)$ be a polynomial over a field $F$ given by $(1.1)$ and $(1.2)$.
Provided $n\gs k$, for any finite subset $A$ of $F$ we have
$$\aligned&|\{f(x_1,\ldots,x_n):\, x_1,\ldots,x_n\in A,\ \t{and}\ x_i\not=x_j\ \t{if}\ i\not=j\}|
\\\gs&\min\l\{p(F)-[\![n=2\ \&\ a_1=-a_2]\!], \f{n(|A|-n)-\{n\}_k\{|A|-n\}_k}k+1\r\}.\endaligned\tag1.5$$
\endproclaim

Sun [S08] noted that this conjecture in the case $n=2$ follows from [PS, Corollary 3], and proved (1.5) with
the lower bound replaced by $\min\{p(F),\,|A|-n+1\}$.

 In this paper we establish a similar version of (1.4) for the case $n\gs k$
 under the condition $a_1=\cdots=a_n$. It implies Conjecture 1.1 in the case
 $a_1=\cdots=a_n$.

 Here is our first result.
\proclaim{Theorem 1.2} Let $A_1, \ldots, A_n$ be finite  subsets
of a field $F$ with $|A_i|\gs i$ for $i=1,\ldots,n$. Let
$$f(x_1,\ldots,x_n)=x_1^k+\cdots+x_n^k+g(x_1,\ldots,x_n)\in
F[x_1,\ldots,x_n]\tag1.6$$ with $\deg g<k\ls n$. Then
$$\aligned&|\{f(x_1,\ldots,x_n):\, x_1\in A_1,\ldots,x_n\in A_n, \t{and}\ x_i\not=x_j\text{ if
}i\not=j\}|
\\&\qquad\gs\min\{p(F),\, q_1+\cdots+q_n+1\}
\endaligned\tag1.7$$
where $$q_i=\min\Sb i\ls j\ls n\\j\eq i\ (\mo\
k)\endSb\l\lfloor\f{|A_j|-j}k\r\rfloor\quad\t{for}\
i=1,\ldots,n.\tag1.8$$
\endproclaim

\Remark\ 1.1. If $k\gs n$ then $q_i=\lfloor(|A_i|-i)/k\rfloor$ for $i=1,\ldots,n$. So
Theorem 1.2 is a complement to Theorem 1.1(ii). In the case $k=1$, Theorem 1.2 yields the main result
of [ANR].

\medskip

Theorem 1.2, together with Theorem 1.1(ii), implies the following extension of the Erd\H os-Heilbronn conjecture.

\proclaim{Theorem 1.3} Let $A_1, \ldots, A_n$ be finite  subsets
of a field $F$ with $|A_1|=\cdots=|A_n|=m\gs n$, and let
$f(x_1,\ldots,x_n)$ be given by $(1.6)$ with $\deg g<k$. Then
$$\aligned&|\{f(x_1,\ldots,x_n):\, x_1\in A_1,\ldots,x_n\in A_n,\ \t{and}\ x_i\not=x_j\text{ if}\ i\not=j\}|
\\\gs&\min\bg\{p(F),\, \f{n(m-n)-\{n\}_k\{m-n\}_k}k+\{n\}_k[\![\{m\}_k<\{n\}_k]\!]+1\bg\}.
\endaligned\tag1.9$$
\endproclaim

\Remark\ 1.2. If  $n$ or $m-n$ is divisible by $k$, then the lower bound in (1.9) becomes
$\min\{p(F),n(m-n)/k+1\}$. In the case $k=1$ and $A_1=\cdots=A_n$, Theorem 1.3 yields
the Dias da Silva-Hamidoune extension (cf. [DH]) of the Erd\H os-Heilbronn conjecture.

In the next section we are going to present an auxiliary theorem.
Theorems 1.2 and 1.3 will be proved in Section 3.

\heading{2. An Auxiliary Theorem}\endheading

For a polynomial $P(x_1,\ldots,x_n)$ over a field, by
$[x_1^{k_1}\cdots x_n^{k_n}]P(x_1,\ldots,x_n)$ we mean the
coefficient of the monomial $x_1^{k_1}\cdots x_n^{k_n}$ in
$P(x_1,\ldots,x_n)$.

In this section we prove the following auxiliary result.

\proclaim{Theorem 2.1} Let $q_1,\ldots, q_n\in\N=\{0,1,2,\ldots\}$
and $k\in\{1,\ldots,n\}$. Then
$$\aligned&\bigg[\prod_{j=1}^nx_j^{kq_j+j-1}\bigg](x_1^k+\cdots+x_n^k)^{N}\prod_{1\ls
i<j\ls n}(x_j-x_i)\\=&N!\prod_{s=1}^k \f{\prod_{0\ls i<j\ls
\lfloor(n-s)/k\rfloor}(q_{jk+s}+j-(q_{ik+s}+i))}
{\prod_{j=0}^{\lfloor(n-s)/k\rfloor}(q_{jk+s}+j)!},
\endaligned\tag2.1$$
where $N=q_1+\cdots+q_n$.
\endproclaim

To prove Theorem 2.1, we need a lemma.

\proclaim{Lemma 2.1} Let $\sigma$ be a
permutation of a finite nonempty set $X$. Suppose that $A$ is a subset of $X$ with
$\sigma(A)=A$. Then
$$\varepsilon(\sigma)=\varepsilon(\sigma|_A)\varepsilon(\sigma|_{X\setminus A}),\tag2.2$$
where $\varepsilon(\sigma)$ stands for the sign of $\sigma$ and
$\sigma|_A$ denotes the restriction of $\sigma$ on $A$.
\endproclaim
\Proof.  Write $\sigma=\tau_1\tau_2\cdots\tau_k$,
where $\tau_1, \tau_2, \ldots, \tau_k$ are disjoint cycles.
As $\sigma(A)=A$, for each $i=1,\ldots,k$, either
all elements in the cycle $\tau_i$ lie in $A$, or
none of the elements in the cycle $\tau_i$ belongs to $A$. Set
$$I=\{1\leqslant i\leqslant k:\,\t{all the elements in the cycle}\ \tau_i\ \t{lie in}\ A\}$$
and
$$\bar{I}=\{1\leqslant i\leqslant k:\,\t{all the elements in the cycle}\ \tau_i\ \t{lie in}\ X\sm A\}$$
Then
$$I\cup \bar I=\{1,\ldots,k\},\ \sigma|_A=\prod_{i\in
I}\tau_i|_A\ \text{and}\ \sigma|_{X\setminus
A}=\prod_{i\in\bar{I}}\tau_i|_{X\sm A}.$$
Therefore
$$\varepsilon(\sigma)=\prod_{i=1}^k\ve(\tau_i)
=\prod_{i\in I}\ve(\tau_i)\times\prod_{j\in \bar I}\ve(\tau_j)
=\varepsilon(\sigma|_A)\varepsilon(\sigma|_{X\setminus A}).$$
This completes the proof.
\qed

\medskip

For a finite nonempty set $X$, we let $S(X)$ denote the symmetric group
of all permutations of $X$. If $|X|=n$, then the group $S(X)$ is isomorphic to
the symmetry group $S_n=S(\{1,\ldots,n\})$. Recall that the determinant of a matrix $[a_{i,j}]_{1\ls i,j\ls n}$
over a field is defined as
$$\det[a_{i,j}]_{1\ls i,j\ls n}=\sum_{\sigma\in S_n}\ve(\sigma)\prod_{i=1}^na_{i,\sigma(i)}
=\sum_{\sigma\in S_n}\ve(\sigma)\prod_{j=1}^na_{\sigma(j),j}.$$

\medskip

\noindent{\it Proof of Theorem 2.1}. By the multinomial theorem,
$$(x_1^k+\cdots+x_n^k)^N=\sum\Sb i_1,\ldots,i_n\in\N\\i_1+\cdots+i_n=N\endSb\frac{N!}{i_1!\cdots
i_n!}x_1^{ki_1}\cdots x_n^{ki_n}.$$
In view of linear algebra,
$$\sum_{\sigma\in S_n}\ve(\sigma)\prod_{j=1}^nx_j^{\sigma(j)-1}
=\det[x_j^{i-1}]_{1\ls i,j\ls n}
=\prod_{1\ls i<j\ls n}(x_j-x_i)\ (\t{Vandermonde)}.$$
Thus
$$ \align
&\bigg[\prod_{j=1}^nx_j^{kq_j+j-1}\bigg](x_1^k+\cdots+x_n^k)^N\prod_{1\ls
i<j\leqslant n}(x_j-x_i)
\\=&\bigg[\prod_{j=1}^nx_j^{kq_j+j-1}\bigg]N!\sum\Sb i_1,\ldots,i_n\in\N\\i_1+\cdots+i_n=N\endSb\sum_{\Sp\sigma\in S_n\\
k\mid \sigma(j)-j\\
\t{for}\ j=1,\ldots,n\endSp}\varepsilon(\sigma)
\prod_{j=1}^n\frac{x_j^{ki_j+\sigma(j)-1}}{i_j!}
\\=&\bigg[\prod_{j=1}^nx_j^{kq_j+j-1}\bigg]N!\sum\Sb i_1,\ldots,i_n\in\N\\i_1+\cdots+i_n=N\endSb
\sum_{\Sp\sigma\in S_n\\
\sigma(X_s)=X_s\\
\t{for}\ s=1,\ldots,k\endSp}\varepsilon(\sigma)
\prod_{j=1}^n\frac{x_j^{ki_j+\sigma(j)-1}}{i_j!},
\endalign
$$
where
$$X_s=\{1\leqslant j\leqslant n:\,j\eq s\ (\mo\ k)\}.$$

 Set $n_s=|X_s|$. Then
$$n_s=|\{q\in\N:\, s+kq\ls n\}|=\l\lfloor\f{n-s}k\r\rfloor+1.$$
If $\sigma\in S_n$ and $\sigma(X_s)=X_s$ for all $s=1,\ldots,k$, then
$$\varepsilon(\sigma)=\ve(\sigma|_{X_1})\ve(\sigma|_{X_2\cup\cdots\cup X_k})=\cdots=
\ve(\sigma|_{X_1})\cdots\ve(\sigma|_{X_k})$$
by Lemma 2.1.

 Let $(x)_0=1$ and
$(x)_i=\prod_{r=0}^{i-1}(x-r)$ for $i=1,2,3,\ldots$.
By the above,
$$\align
&\bigg[\prod_{j=1}^nx_j^{kq_j+j-1}\bigg]\f{(x_1^k+\cdots+x_n^k)^N}{N!}\prod_{1\leqslant
i<j\leqslant n}(x_j-x_i)\\
=&\bigg[\prod_{j=1}^nx_j^{kq_j+j-1}\bigg]\sum\Sb
i_1,\ldots,i_n\in\N\\i_1+\cdots+i_n=N\endSb
\prod_{s=1}^{k}\bigg(\sum_{\sigma\in
S(X_s)}\varepsilon(\sigma)\prod_{j\in
X_s}\frac{x_j^{ki_j+\sigma(j)-1}}{i_j!}
\bigg)\\
=&\prod_{s=1}^{k}\bigg(\sum\Sb\sigma\in
S(X_s)\\q_j+(j-\sigma(j))/k\gs0\\\t{for}\ j\in X_s\endSb\varepsilon(\sigma)\prod_{j\in
X_s}\frac{1}{(q_j+(j-\sigma(j))/k)!}\bigg)
\\=&\prod_{s=1}^{k}\bigg(\sum_{\sigma\in
S_{n_s}}\varepsilon(\sigma)\prod_{j=1}^{n_s}\frac{(q_{(j-1)k+s}+j-1)_{\sigma(j)-1}}{(q_{(j-1)k+s}+j-1)!}\bigg)
\\=&\prod_{s=1}^{k}\f{\det[(q_{(j-1)k+s}+j-1)_{i-1}]_{1\ls i,j\ls n_s}}{\prod_{j=1}^{n_s}(q_{(j-1)k+s}+j-1)!}
=\prod_{s=1}^{k}\f{\det[(q_{jk+s}+j)_i]_{0\ls i,j\ls n_s-1}}{\prod_{j=0}^{n_s-1}(q_{jk+s}+j)!}.
\endalign$$
It is well known that
$$y^{i}=(y)_{i}+\sum_{0\ls r<i}S(i,r)(y)_{r}\quad\t{for}\ i=0,1,2,\ldots,$$
where $S(i,r)\ (0\ls r<i)$ are Stirling numbers of the second kind.
So
$$\align &\det[(q_{jk+s}+j)_{i}]_{0\ls i,j<n_s}=\det[(q_{jk+s}+j)^{i}]_{0\ls i,j<n_s}
\\=&\prod_{0\ls i<j<n_s}(q_{jk+s}+j-(q_{ik+s}+i))\ \ (\t{Vandermonde}),
\endalign$$
and hence (2.1) follows. \qed

\heading{3. Proofs of Theorems 1.2 and 1.3}\endheading

Let us recall the following powerful tool.

\medskip

\noindent{\bf Combinatorial Nullstellensatz} (Alon \cite{A}). {\it
Let $A_1,\cs,A_n$ be finite subsets of a field $F$, and let
$P(x_1,\cs,x_n)\in F[x_1,\cs,x_n]$. Suppose that $\deg
P=k_1+\cdots+k_n$ where $0\ls k_i<|A_i|$ for $i=1,\ldots,n$.
If
$$[x_1^{k_1}\cdots x_n^{k_n}]P(x_1,\ldots,x_n)\not=0,$$
then
$P(x_1,\ldots,x_n)\not=0$ for some $x_1\in A_1,\ldots,x_n\in A_n$.}

\medskip

\noindent{\it Proof of Theorem 1.2}.
Let $m$ be the least nonnegative integer not exceeding $n$
such that $\sum_{m<i\ls n}q_i<p(F)$.
For each $m<i\ls n$ let $A_i'$ be a subset of $A_i$ with cardinality $kq_i+i\ls|A_i|$.
In the case $m>0$, $p=p(F)$ is a prime and we let $A_{m}'$ be a subset of $A_{m}$
with
$$|A_{m}'|=k\(p-1-\sum_{m<i\ls n}q_i\)+m<kq_{m}+m\ls|A_{m}|.$$
If $0<i<m$ then we let $A_{i}'\se A_i$ with $|A_i|=i$.
Clearly $q_i'=(|A_i'|-i)/k\ls q_i$.
Whether $m=0$ or not, we have
$\sum_{i=1}^n(|A_i'|-i)=k\sum_{i=1}^nq_i'=k(N-1),$ where
$$N=\min\{p(F),\,q_1+\cdots+q_n+1\}.$$

Let $s\in\{1,\ldots,k\}$. For any  $0<i<n_s=\lfloor(n-s)/k\rfloor+1$ we have
$$q_{(i-1)k+s}=\min\Sb (i-1)k+s\ls j\ls n\\j\eq s\,(\mo\ k)\endSb\l\lfloor\f{|A_j|-j}k\r\rfloor
\ls\min\Sb ik+s\ls j\ls n\\j\eq s\,(\mo\ k)\endSb\l\lfloor\f{|A_j|-j}k\r\rfloor=q_{ik+s}$$
and hence $q_{(i-1)k+s}'\ls q_{ik+s}'$. (If $(i-1)k+s=m$ then $q_{(i-1)k+s}'\ls q_{(i-1)k+s}\ls q_{ik+s}=q_{ik+s}'$.)
So
$$0\ls q_s'<q_{k+s}'+1<q_{2k+s}'+2<\cdots<q_{(n_s-1)k+s}'+n_s-1.$$

Define
$$P(x_1,\ldots,x_n)=(x_1^k+\cdots+x_n^k)^{N-1}\prod_{1\ls i<j\ls n}(x_j-x_i)\in F[x_1,\ldots,x_n].$$
In light of Theorem 2.1,
$$\bg[\prod_{j=1}^nx_j^{kq_j'+j-1}\bg]P(x_1,\ldots,x_n)=he,$$
where $e$ is the identity of the field $F$, and
$$h=(N-1)!\bigg/\prod_{s=1}^k\prod_{j=0}^{n_s-1}
\prod\Sb 0\ls r<q_{jk+s}'+j\\r\not\in\{q_{ik+s}'+i:\,0\ls i<j\}\endSb(q_{jk+s}'+j-r)$$
is an integer dividing $(N-1)!$. Since $p(F)>N-1$, we have
$he\not=0$.

Set
$$C=\{f(x_1,\ldots,x_n):\, x_1\in A_1',\ldots,x_n\in A_n',\,\t{and}\ x_i\not=x_j\ \t{if}\ i\not=j\}.$$
Suppose that $|C|\ls N-1$ and let $Q(x_1,\ldots,x_n)$ denote the polynomial
$$f(x_1,\ldots,x_n)^{N-1-|C|}\prod_{c\in C}(f(x_1,\ldots,x_n)-c)
\times\prod_{1\ls i<j\ls n}(x_j-x_i).$$
Then
$$\deg Q=k(N-1)+\bi n2=\sum_{i=1}^n(|A_i'|-1)=\sum_{i=1}^n(kq_i'+i-1)$$
and
$$\l[x_1^{|A_1'|-1}\cdots x_n^{|A_n'|-1}\r]Q(x_1,\ldots,x_n)
=\bg[\prod_{j=1}^nx_j^{kq_j'+j-1}\bg]P(x_1,\ldots,x_n)\not=0.$$
In light of the Combinatorial Nullstellensatz, there are $x_1\in A_1',\ldots,x_n\in A_n'$ such that
$Q(x_1,\ldots,x_n)\not=0$. This contradicts the fact $f(x_1,\ldots,x_n)\in C$.

 By the above, we have
 $$|\{f(x_1,\ldots,x_n):\, x_1\in A_1,\ldots,x_n\in A_n,\,\t{and}\ x_i\not=x_j\ \t{if}\ i\not=j\}|
 \gs|C|\gs N.$$
 This concludes the proof. \qed

\medskip
\noindent{\it Proof of Theorem 1.3}. Write $n=kq_0+n_0$ with $q_0\in\N$ and $1\ls n_0\ls k$.
Then
$$\align&\sum_{i=1}^n\min\Sb i\ls j\ls n\\j\eq i\,(\mo\ k)\endSb\l\lfloor\f{|A_j|-j}k\r\rfloor
\\=&\sum_{r=1}^{k}\ \sum_{0\ls q\ls\lfloor(n-r)/k\rfloor}\min\Sb kq+r\ls j\ls n\\j\eq r\,(\mo\ k)\endSb
\l\lfloor\f{m-j}k\r\rfloor
\\=&\sum_{r=1}^k\ \sum_{0\ls q\ls\lfloor(n-r)/k\rfloor}\l\lfloor\f{m-r-k\lfloor(n-r)/k\rfloor}k\r\rfloor
\\=&\sum_{r=1}^k\(\l\lfloor\f{n-r}k\r\rfloor+1\)\(\l\lfloor\f{m-r}k\r\rfloor-\l\lfloor\f{n-r}k\r\rfloor\)
\\=&\sum_{r=1}^{n_0}(q_0+1)\(\l\lfloor\f{m-r}k\r\rfloor-q_0\)
+\sum_{n_0<r\ls k}q_0\(\l\lfloor\f{m-r}k\r\rfloor-q_0+1\)
\\=&q_0\sum_{r=1}^k\l\lfloor\f{m-r}k\r\rfloor+\sum_{r=1}^{n_0}\l\lfloor\f{m-r}k\r\rfloor
-q_0((q_0+1)n_0+(q_0-1)(k-n_0)).
\endalign$$
Observe that
$$\sum_{r=1}^k\l\lfloor\f{m-r}k\r\rfloor
=\sum_{r=1}^k\(\f{m-r}k-\f{\{m-r\}_k}k\)=m-\sum_{r=1}^k\f rk-\sum_{s=0}^{k-1}\f sk=m-k.$$
So we have
$$\align&\sum_{i=1}^n\min\Sb i\ls j\ls n\\j\eq i\,(\mo\ k)\endSb\l\lfloor\f{|A_j|-j}k\r\rfloor
\\=&\sum_{r=1}^{n_0}\l\lfloor\f{m-r}k\r\rfloor+q_0(m-k)-q_0(2n_0+k(q_0-1))
\\=&\sum_{r=1}^{n_0}\(\l\lfloor\f{m-r}k\r\rfloor-q_0\)+q_0(m-n).
\endalign$$

To simplify the last sum, we now subtract a term which will be added later.
Clearly
$$\align&\sum_{r=1}^{n_0}\(\l\lfloor\f{m-r}k\r\rfloor-q_0\)-n_0\l\lfloor\f{m-n}k\r\rfloor
\\=&\sum_{r=1}^{n_0}\(\l\lfloor\f{m-n+n_0-r}k\r\rfloor-\l\lfloor\f{m-n}k\r\rfloor\)
\\=&\sum_{s=0}^{n_0-1}\l\lfloor\f{\{m-n\}_k+s}k\r\rfloor.
\endalign$$
If $\{m\}_k\gs n_0$, then
$$\sum_{s=0}^{n_0-1}\l\lfloor\f{\{m-n\}_k+s}k\r\rfloor
=\sum_{s=0}^{n_0-1}\l\lfloor\f{\{m\}_k-\{n\}_k+s}k\r\rfloor=0.$$
If $\{m\}_k<n_0$, then
$$\align &\sum_{s=0}^{n_0-1}\l\lfloor\f{\{m-n\}_k+s}k\r\rfloor
=\sum_{s=0}^{n_0-1}\l\lfloor\f{\{m\}_k-n_0+k+s}k\r\rfloor
\\=&|\{s\in\{0,\ldots,n_0-1\}:\ s\gs n_0-\{m\}_k\}|=\{m\}_k.
\endalign$$
Therefore
$$\align&\sum_{i=1}^n\min\Sb i\ls j\ls n\\j\eq i\,(\mo\ k)\endSb\l\lfloor\f{|A_j|-j}k\r\rfloor
\\=&q_0(m-n)+n_0\l\lfloor\f{m-n}k\r\rfloor+\{m\}_k[\![\{m\}_k<n_0]\!]
\\=&(m-n)\l\lfloor\f nk\r\rfloor+\{n\}_k\l\lfloor\f{m-n}k\r\rfloor+\{m\}_k[\![\{m\}_k<\{n\}_k]\!]
\\=&\f{n(m-n)}k-\f{\{n\}_k\{m-n\}_k}k+\{m\}_k[\![\{m\}_k<\{n\}_k]\!].
\endalign$$

 In view of the above, by applying Theorem 1.2 for $k\ls n$ and Theorem 1.1(ii) for $k\gs n$,
we immediately get the desired (1.9). \qed


\medskip

\widestnumber \key{ANR}

\Refs

\ref\key A \by N. Alon\paper Combinatorial Nullstellensatz\jour
Combin. Probab. Comput.\vol8\yr1999\pages7--29\endref

\ref\key ANR\by N. Alon, M. B. Nathanson and I. Z. Ruzsa \paper
The polynomial method and restricted sums of congruence classes
\jour J. Number Theory\vol56\yr1996\pages404--417\endref

\ref\key BW\by P. Balister and J. P. Wheeler\paper
 The Erd\H os-Heilbronn conjecture for finite groups
\jour Acta Arith.\pages to appear\endref

\ref\key DH\by J. A. Dias da Silva and Y. O. Hamidoune\paper
Cyclic spaces for Grassmann derivatives and additive theory\jour
Bull. London Math. Soc.\vol 26\yr1994\pages140--146\endref

\ref\key EH\by P. Erd\H os and H. Heilbronn \paper On the addition
of residue classes modulo $p$ \jour Acta
Arith.\vol9\yr1964\pages149--159\endref

\ref\key HS\by Q. H. Hou and Z. W. Sun\paper Restricted sums in a
field \jour Acta Arith.\vol 102\yr 2002\pages 239--249\endref

\ref\key N\by M. B. Nathanson\book Additive Number Theory: Inverse
Problems and the Geometry of Sumsets \publ Graduate Texts in
Math., Vol. 165, Springer, New York\yr 1996\endref

\ref\key PS\by H. Pan and Z. W. Sun\paper A lower bound for
$|\{a+b\colon a\in A,\ b\in B,\ P(a,b)\not=0\}|$ \jour J. Combin.
Theory Ser. A\vol 100\yr 2002\pages 387--393\endref

\ref\key S03\by Z. W. Sun\paper On Snevily's conjecture and
restricted sumsets\jour J. Combin. Theory Ser. A \vol 103\yr 2003
\pages 291--304\endref

\ref\key S08\by Z. W. Sun\paper On value sets of polynomials over a
field \jour{\it Finite Fields Appl.}\vol 14\yr 2008\pages
470--481\endref

\ref\key S08a\by Z. W. Sun\paper An additive theorem and restricted
sumsets \jour{\it Math. Res. Lett.}\vol 15\yr 2008\pages
1263--1276\endref

\ref\key SY\by Z. W. Sun and Y. N. Yeh\paper On various restricted sumsets
\jour J. Number Theory\vol 114\yr 2005\pages 209--220\endref

\ref\key TV\by T. Tao and V. H. Vu\book Additive Combinatorics
\publ Cambridge Univ. Press, Cambridge, 2006\endref

\endRefs
\enddocument